\documentclass[]{siamart0516}
\usepackage[utf8]{inputenc}

\usepackage{listings}
\usepackage{changes}
\definechangesauthor[name=Stephan, color=blue!60!black]{STS}

\title{A Shape Newton Scheme for Deforming Shells with Application to Capillary Bridges}
\author{Stephan Schmidt\footnotemark[1] \footnotemark[4], Melanie Gräßer\footnotemark[2], Hans-Joachim Schmid\footnotemark[2]}

\usepackage{amsmath, amssymb}
\usepackage{bm} 
\usepackage{tikz}
\usetikzlibrary{positioning, shapes, arrows, calc}

\usepackage{pgfplots}
\pgfplotsset{compat=1.16}
\pgfplotsset{}

\usepackage{graphicx,color}
\usepackage{url}
\usepackage{microtype}
\usepackage{caption}
\usepackage{subcaption}
\usepackage{overpic}

\usepackage[style=numeric,sorting=nyt,giveninits=true,maxbibnames=99,url=false]{biblatex}
\DeclareFieldFormat[article]{titlecase}{\MakeSentenceCase{#1}}

\def\clap#1{\hbox to 0pt{\hss#1\hss}}

\usepackage{mathrsfs}

\usepackage{stmaryrd}

\newcommand{\R}{\mathbb{R}}

\newcommand{\norm}[1]{\left\Vert#1\right\Vert}

\newcommand{\inner}[2]{\left\langle {#1}, {#2}\right\rangle}

\def\bei#1{\vrule width 0.4pt height 14pt depth 9pt
           \lower 8pt \hbox{$ _{\hbox{}\, #1}$}\!}
\renewcommand{\div}[1][]{{\operatorname{div}_{#1} \,}} 
\def \tr{\operatorname{tr}}

\newcommand{\imagefolder}{./Images}

\def \dist {\operatorname{dist}} 
\def \dd {\operatorname{d}} 
\def \dD{\operatorname{D}}

\def \bd {\mathbf{d}}

\def \bF {\mathbf{F}}
\def \bf {\mathbf{f}}
\def \bg {\mathbf{g}}
\def \bH {\mathbf{H}}
\def \bN {\mathbf{N}}
\def \bn {\mathbf{n}}
\def \bu {\mathbf{u}}
\def \bV {\mathbf{V}}
\def \bW {\mathbf{W}}

\def \bmu {\bm{\mu}}

\begin{document}

\definecolor{dkgreen}{rgb}{0,0.6,0}
\definecolor{gray}{rgb}{0.5,0.5,0.5}
\definecolor{mauve}{rgb}{0.58,0,0.82}

\lstset{frame=tb,
  language=Python,
  aboveskip=3mm,
  belowskip=3mm,
  showstringspaces=false,
  columns=flexible,
  basicstyle={\small\ttfamily},
  numbers=none,
  numberstyle=\tiny\color{gray},
  keywordstyle=\color{blue},
  commentstyle=\color{dkgreen},
  stringstyle=\color{mauve},
  breaklines=true,
  breakatwhitespace=true,
  tabsize=3
}
\lstMakeShortInline[columns=fixed]|
 \tikzstyle{block} = [rectangle, draw, minimum width=8em, align=center, rounded corners, minimum height=2em, scale=0.9]
\tikzstyle{line} = [draw, -latex']
\tikzstyle{deriv} = [draw, -latex', color=red]
\tikzstyle{hderiv} = [draw, -latex', scale=0.9]

\maketitle

\renewcommand{\thefootnote}{\fnsymbol{footnote}}
\footnotetext[1]{Department of Mathematics, Humboldt University Berlin, 10099 Berlin, Germany} 
\footnotetext[2]{Faculty of Mechanical Engineering, Paderborn University, 33098 Paderborn, Germany} 
\footnotetext[4]{Corresponding author email: ({\tt s.schmidt@hu-berlin.de})} 
\renewcommand{\thefootnote}{\arabic{footnote}}


\begin{abstract}
We present a second order numerical scheme to compute capillary bridges between arbitrary solids by minimizing the total energy of all interfaces. From a theoretical point of view, this approach can be interpreted as the computation of generalized minimal surfaces using a Newton-scheme utilizing the shape Hessian. In particular, we give an explicit representation of the shape Hessian for functionals on shells involving the normal vector without reverting back to a volume formulation or approximating curvature. From an algorithmic perspective, we combine a resolved interface via a triangulated surface for the liquid with a level set description for the constraints stemming from the arbitrary geometry. The actual shape of the capillary bridge is then computed via finite elements provided by the FEniCS environment, minimizing the shape derivative of the total interface energy.
\end{abstract}

\begin{keywords}
Shape Optimization, Shape Newton, Shape Hessian, Shells, Capillary Bridges
\end{keywords}

\begin{AMS}
49M05, 49M15, 49M37, 65D18, 65K10
\end{AMS}

\pagestyle{myheadings}
\thispagestyle{plain}

\section{Introduction}
In process engineering the behavior of particle systems is essential for the design of many processes such as  particle agglomeration and separation, powder handling and granular flow~\cite{butt_2009, forsyth_2002, landi_2011}. A well developed method for the simulation of particle systems is the discrete element method, which calculates the system parameters by solving Newton's equation of motion for each particle. Thus, suitable force models for the particle interaction are indispensable. Besides Van-der-Waals forces and electrostatic forces, capillary forces play a dominant role in wet materials \cite{rumpf_1974}, as liquid bridges form between particles. These bridges add forces to the system due to the pressure difference between the liquid and gas phase, as well as due to surface tension. For very small particles, even condensation out of humid air may lead to capillary bridges, which cause a significant attractive force, typically exceeding other inter-particle forces. The curved phase boundary between the liquid and the surrounding gas leads to a pressure difference over the interface. The mean curvature $H$ and the pressure difference $\Delta p$ are directly connected by the surface tension $\sigma$. This relation is described by the Young-Laplace-Equation
\begin{align}\tag{YLE}\label{eq:YLE}
\Delta p(s) = -\sigma 2H(s),
\end{align} 
where $s$ is the spatial variable. This relationship can be derived by the balance of forces~\cite{Laplace1805, Young_1805} or by minimizing the system energy~\cite{Gauss1830}. In case of negligible gravity, $\Delta p$ becomes spatially independent and a constant mean curvature (CMC) surface describes the equilibrium state.

From an application oriented point of view, the methods can be divided into general methods for arbitrary geometries and two dimensional axisymmetric methods. Most of the literature focuses on axisymmetric bridges between two ideal spheres, a sphere and a plane or other very specific geometries. This also implies a constant contact angle and negligence of gravity~\cite{Schubert1982}. Surfaces of revolution with a constant mean curvature have already been studied in 1841 by Delaunay~\cite{Delaunay1841}, who expressed them in terms of a non-linear ordinary differential equation, which is derived by rolling a conic along a straight line. The theory of CMC surfaces was first applied to capillary bridges by Plateau~\cite{Plateau1864}. He showed that with rising the volume of the liquid, a sequence of different CMC surface types occurs. This sequence was also analyzed by Orr, Scriven and Rivas~\cite{orr_scriven_rivas_1975}, who solved~\eqref{eq:YLE} in two dimensions using elliptic integrals for a sphere and a plate. Later, Rubinstein and Fel~\cite{rubinstein_2014} proofed that the classical Plateau sequence is only valid if the sphere and the plate are touching each other. In case of a non-zero distance, the sequence might change. Moreover, for special parameter combinations, the solution of the elliptic integrals might get complex, which leads to multiple solutions of the system that require further stability analysis. Besides solving~\eqref{eq:YLE} in terms of elliptic integrals, shooting methods have been presented in~\cite{dormann_2015, qiangNian_2017}. Moreover, several methods requiring simplification beyond rotational symmetry, such as the toroidal approximation~\cite{fisher_1926} and the parabolic approximation~\cite{pepin_2000} are described in the literature. Methods not tailor made to exploit the axisymmetric situation can still require or exploit specific settings. For example, molecular dynamic methods~\cite{ko_2010, laube_2017} are restricted to very small bridges with a manageable number of molecules, whereas computation of the bridge as the steady state of a free surface computational fluid dynamics approach requires an accurate representation of the surface~\cite{jettestuen_2013, sun_2016}.

The general idea of energy based approaches is to find the steady state solution of the system by minimizing its surface energy. Treating the formation of capillary bridges via variational formulation and optimization has been discussed in~\cite{HornungMittelmann1990}. Most, if not all of the energy minimization methods first discretize the unknown surface of the capillary bridge via triangles and then seek to minimize the discrete reformulation of the system energy as a function of the finite vertex positions. However, very often the continuous or ``classical'' formulation of~\eqref{eq:YLE} is still utilized in these schemes, resulting in the necessity to somehow find an approximation of the mean curvature, which is not readily available when using triangulated meshes. There are several discretized energy minimization methods for the calculation of CMC surfaces in general geometries, see, e.g.~\cite{chen_2008, dziuk_2006, polthier_2002}. Nevertheless, successful application of any of these methods strongly depends on a high mesh quality including a uniform vertex distribution and minor triangle deformation. Thus, the minimization can also be done on a least square functional which is easier to handle compared to the surface area functional and allows to use an extended centroidal Voronoi tessellation~\cite{pan_2012}. Using this method, the mesh can be optimized simultaneously with the surface area. A modified approach is presented by Renka~\cite{renka_1995, renka_2015}, who derived a non-linear least square system by using a different energy functional. The major advantage of his method is the robustness in the presence of low mesh quality. Another option of solving the capillary bridge problem by energy minimization is the direct evaluation of virtual displacements \cite{iliev_1995}. 
Additionally,  Ardito~et~al.~\cite{ardito_2014} applied the mechanical concept of an ideal stiff thin membrane. Based on the total potential energy, the problem is expressed according to the classical principle of virtual power. The results of axisymmetric bridges agree well with CMC calculations of Kenmotsu~\cite{kenmotsu_1980} while deviations from molecular dynamic simulations can be found in Ko~et~al.~\cite{ko_2010}. Although this method is capable of calculating  asymmetric cases no examples are outlined.

A frequently used open source code following this paradigm is \textit{The Surface Evolver}~\cite{brakke_1992}, which also features a Newton scheme using the Hessian of the optimization problem in terms of vertex coordinates~\cite{brakke_1997}. The surface energy might include surface tension, gravity or other energy forms. In addition, surface constraints, such as a specific geometry or body volume, can be specified. The general implementation of \textit{The Surface Evolver} is a fundamental advantage and consequently it has been used to study different aspects of the capillary bridge problem such as non-symmetric bridges~\cite{bedarkar_2009, broesch_2012, chau_2007,farmer_2015, virozub_2009}, stability~\cite{ataei_2017} and the influence of gravity~\cite{sun_2018}. 

As part of this work, we present a novel method for minimizing~\eqref{eq:YLE} using a shape Newton scheme. To this end, we use shape calculus to derive a variational reformulation of~\eqref{eq:YLE}, which does not involve curvature and is hence better suited for the lower regularity of triangulated surfaces. Indeed, in comparison to most works using discrete surfaces to minimize the total energy by moving vertices, we derive an explicit, curvature-free expression of the first and second order derivatives of the energy with respect to shape perturbations. Furthermore, we seek expressions on shell meshes only, which still adhere to the low regularity of a constant outer normal per facet of the mesh, that is the concentration of curvature to edges. This leads to several novel expressions for the shape Hessian of surface integrals and the second material derivative of the outer normal, which are also of considerable interest for general problems beyond capillary bridges. Indeed, we provide a general expression of the shape Hessian for boundary integrals involving the outer normal. The shape Hessian for the much simpler iso-perimeter problem has previously been considered in~\cite{SAD_2018}, albeit using either curvature or volume integrals.

The resulting variational, curvature-free alternative representation of~\eqref{eq:YLE} and the corresponding shape Hessian can then be solved by any kind of finite element software offering continuous Lagrange elements on shells. To this end, we use the FEniCS environment~\cite{Fenics_book}, which has full support for shell meshes~\cite{FEniCSManifold}. In particular, FEniCS also offers the automatic computation of discrete shape derivatives for vertex motions, provided they do not involve cross-terms such as in the off-diagonal blocks of the shape Hessian~\cite{FenicsShapeAD}. We use this functionality to confirm that the ``optimize-then-discretize'' approach presented here commutes with the ``discretize-then-optimize'' paradigm for all first and second order partial derivatives with respect to the shape only, irrespective of mesh resolution.

This paper is structured as follows. Section~\ref{sec:Problem} summarizes shape calculus simultaneously to deriving a curvature free variational reformulation of the Young-Laplace equation for triangulated meshes, which corresponds to the shape derivative of the energy of the system without additional geometric constraints such as the touch constraint of the liquid to the solid body surfaces. Section~\ref{sec:ShapeHessian} then focuses on the shape Hessian of the problem and we derive the second order shape derivative of the Lagrangian of the energy minimization problem. Special attention is again given to the fact that no curvature and no volume integrals arise. This necessitates the computation of second order material derivatives of the normal, but consequently also leads to several novel expressions of the shape Hessian for boundary integrals, which are also of interest in general problems beyond capillary bridges. Section~\ref{sec:Implementation} then discusses our numerical implementation in detail, in particular with respect to possible rank deficiencies in the shape Hessian. Finally, we test our numerical implementation by revisiting several standard test cases for capillary bridges from the literature in Section~\ref{sec:NumericalResults}.

\section{Problem Formulation}\label{sec:Problem}
\subsection{Capillary Bridges as Generalized Minimal Surfaces}

\begin{figure}[h!]
\begin{center}
  \begin{overpic}[width=0.7\textwidth]{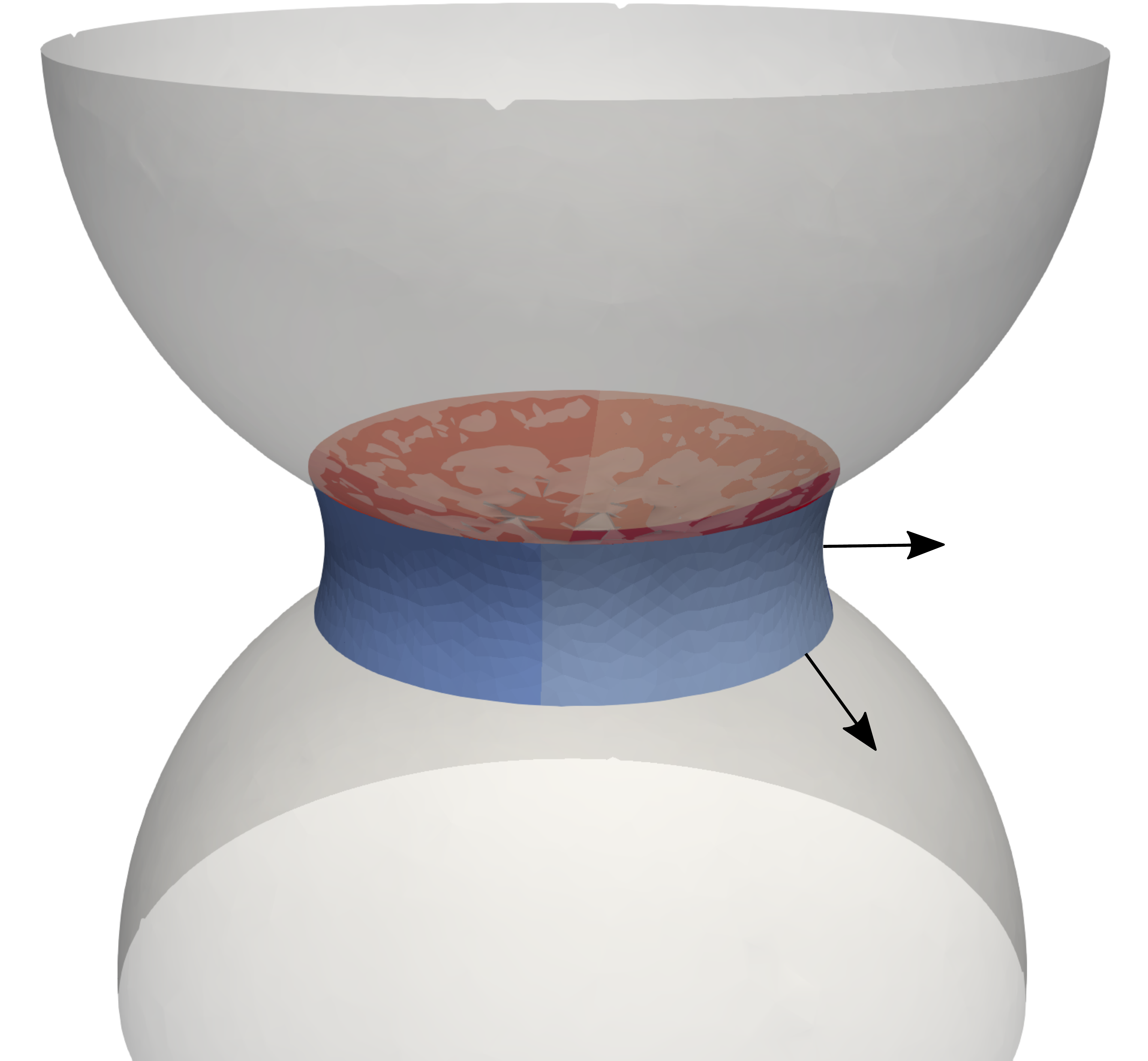}
      \put (20,10){$\Gamma_{\text{S}_1}$}
      \put (30,75){$\Gamma_{\text{S}_2}$}
      \put (83,44){$\bn$}
      \put (77,26){$\bmu_\text{LA}$}
      \put (19,43){$\Gamma_\text{LA}$}
      \put (40,50){$\Gamma_{\text{LS}_2}$}
    \end{overpic}
\end{center}
\caption{The geometrical setup.}\label{fig:setup}
\end{figure}
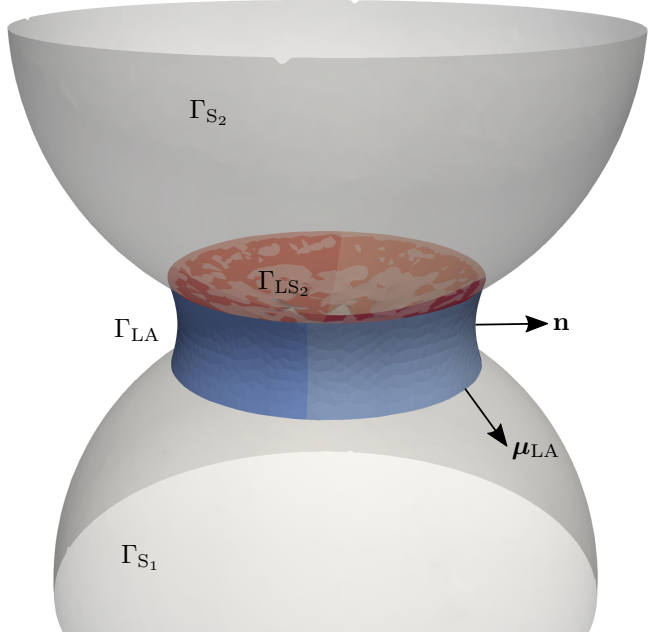

We assume that $\Gamma = \partial \Omega$ is the surface of some volume of fluid $\Omega$. This capillary surface $\Gamma$ can be disjointly decomposed into $\Gamma_{\text{LA}}$, the liquid air interface, and $\Gamma_{\text{LS}_i}$, the interface between the liquid and solid body $i$. The resulting interface lines, or solid-gas-liquid triple phase lines are subsumed as boundaries $\partial \Gamma$. The geometrical setting is shown in Figure~\ref{fig:setup}.
In particular, we use $\bn$ to denote the outer normal of the volume $\Omega$ and $\bmu$ to denote co-normals of $\Gamma$, i.e., unit vectors that are orthogonal to $\partial \Gamma$ and to $\bn$. In order to reduce the computational workload and to avoid additional rank deficits of the Hessian, we are in particular interested in problem formulations that are tailor made for shell meshes, i.e., grids of topological dimension two and geometrical dimension three, and which do not require volume integrals. Hence, the surface area of each interface and the total volume of fluid are then given by the following surface integrals
\begin{align}\label{eq:S}
  S(\Gamma) := \int\limits_{\Gamma} 1 \dd s, \quad F(\Gamma) := \int\limits_{\Omega} 1 \ \dd x = \int\limits_{\Gamma} \frac{1}{3}\inner{s}{\bn}\ \dd s.
\end{align}
Furthermore, the gravitational force acting on the liquid is given by
\begin{align}\label{eq:G}
  G(\Gamma) := \int\limits_{\Omega} \inner{\bg}{x} \dd x = \int\limits_{\Omega} \frac{1}{2}\div((g_ix_i^2)_{i=1}^3) \dd x = \int\limits_{\Gamma} \frac{1}{2} \sum_{i=1}^3g_is_i^2n_i \dd s = \int\limits_{\Gamma}\inner{\tilde{\bg}(s)}{\bn} \dd s,
\end{align}
where $\bg$ is the gravitational vector and $\tilde{\bg}(s) := \frac{1}{2}(g_1s_1^2, g_2s_2^2, g_3s_3^2)^T$.

The formation of the capillary surface follows the minimization of the total interface energy, leading to the problem
\begin{align}\label{eq:main}
  \min_{\Gamma}&\quad S(\Gamma_{\text{LA}}) - \left(\sum_{i=1}^k \beta_i S(\Gamma_{\text{LS}_i}) \right) + b G(\Gamma) \\
  \text{s.t.}&\nonumber\\
  &\nonumber\begin{aligned}
  F(\Gamma) &= V_0\\
  \Gamma_{\text{LS}_i} &\subset \Gamma_{{\text{S}}_i} \quad \forall i=1,...,k,
  \end{aligned}
\end{align}
where, $\beta_i > 0$ is the non-dimensional relative adhesion coefficient of solid $i$ and $b\norm{\bg}\geq 0$ is the Bond number, the non-dimensional gravitational influence. Finally, $\Gamma_{\text{S}_i}$ denotes the surface of solid body $i$.

The interface condition $\Gamma_{\text{LS}_i} \subset \Gamma_{{\text{S}}_i}$ is not directly numerically tractable and can also not readily be addressed via shape calculus.
Suppose the discretization of the subdomain $\Gamma_{\text{LS}_i}$ consists of the vertices $s_{i,j}$, $j=1,...,k_i$. We then replace the analytic subset condition with the constraint of minimal distance for each vertex $s_{i,j}$, namely
\begin{align}\label{eq:DistanceConstraints}
  c(s_{i,j}, \Gamma_{\text{S}_i}) := \dist(s_{i,j}, \Gamma_{\text{S}_i}) = 0, \text{ for all }s_{i,j} \text{ vertex of }\Gamma_{\text{LS}_i},
\end{align}
where $\dist(., \Gamma_{\text{S}_i}): \R^3 \rightarrow \R$ is some signed distance measure from any point in $\R^3$ to $\Gamma_{\text{S}_i}$, the surface of the solid body $i$. We follow the convention of the distance being negative if the point $s$ is in the interior of $\Gamma_{\text{S}_i}$.

\subsection{Shape Calculus}
To find the surface $\Gamma$ solving~\eqref{eq:main}, we employ techniques from shape optimization, in particular shape calculus. Following the approach of perturbation of identity, a deformed surface $\Gamma_\varepsilon$ is defined via
\begin{align*}
  \Gamma_{\varepsilon} := \{s + \varepsilon \bV(s) : s \in \Gamma\},
\end{align*}
where $\bV: \R^3 \rightarrow \R^3$ is a suitably smooth vector field and $s_\varepsilon:=s + \varepsilon \bV(s)$. Following~\cite{DelfourZolesio2, ZolesioSokolowski}, the Eulerian-Semi Derivative of a shape functional
\begin{align*}
  dJ(\Gamma)[\bV] := \frac{\dd}{\dd\varepsilon}\bei{\varepsilon=0^+}J(\Gamma_\varepsilon) =  \frac{\dd}{\dd\varepsilon}\bei{\varepsilon=0^+}\int\limits_{\Gamma_\varepsilon} f(\varepsilon, s_\varepsilon) \ \dd s_\varepsilon
\end{align*}
is then given by
\begin{align}
  dJ(\Gamma)[\bV] &= \int\limits_{\Gamma} f\div[\Gamma] \bV + df[\bV] \ \dd s \label{eq:SurfDerive1} \\
  &= \int\limits_{\Gamma} \div[\Gamma] (f\bV) + \inner{\bV}{\bn} \inner{\nabla f}{\bn} + f'[\bV] \ \dd s \label{eq:SurfDerive2}\\
  &= \int\limits_{\Gamma} \inner{\bV}{\bn} \left[\langle \nabla f, \bn\rangle + \kappa f\right] + f'[\bV] \ \dd s + \int\limits_{\partial \Gamma} \langle \bV, \bmu\rangle f \ \dd \ell,\label{eq:SurfDerive3}
\end{align}
provided the shape is of sufficient regularity, with~\eqref{eq:SurfDerive1} requiring the least. Here, $\div[\Gamma]$ is the divergence in the tangent space of $\Gamma$ and $\kappa := \div[\Gamma] \bn$ is the additive curvature, with $\bn$ being the outer normal to $\Omega$ and $\bmu$ being the normal to $\partial \Gamma$, which also fulfills $\langle \bmu, \bn\rangle = 0$ and is also called the co-normal. Finally, $df[\bV](s) := \frac{d}{d\varepsilon}_{\vert_{\varepsilon=0}} f(\varepsilon, s_\varepsilon)$ is the material derivative and $f'[\bV](s) := \frac{\partial}{\partial\varepsilon}_{\vert_{\varepsilon=0}} f(\varepsilon, s)$ is the local derivative. Due to the chain rule, the identity
\begin{align*}
df[\bV] = f' + \dD f\bV
\end{align*}
holds, if both local and material derivative exist. Here, $\dD$ is the classical spatial Jacobian. Hence, for the spatial coordinates $x$ themselves, one arrives at
\begin{align*}
  dx[\bV] = \dD x\bV = \bV.
\end{align*}
It is worth mentioning that, contrary to local derivatives, the material derivative and spatial differentiation do not commute, meaning
\begin{align}\label{eq:MaterialSpatial}
  \dD df[\bV] = \dD(f'[\bV] + \dD f \bV) = (\dD f)'[\bV] + \dD^2f\bV + \dD f\dD\bV = d(\dD f)[\bV] + \dD f\dD \bV.
\end{align}
For the outer normal, one has in particular~\cite{PHDSchmidtShape, ZolesioSokolowski}
\begin{align}\label{eq:dn}
  d\bn[\bV] = -(\dD_\Gamma\bV)^T\bn,
\end{align}
where $\dD_\Gamma \bV = \dD\bV - \dD\bV\bn\bn^T$ is the tangential Jacobian of $\bV$ with respect to the spatial coordinate.
Hence,~\eqref{eq:SurfDerive2} is essentially the chain rule relating the material derivative to the local derivative, where the tangential part of the convective part has been merged with the divergence. Furthermore, Equation$~\eqref{eq:SurfDerive3}$ is created by integration by parts on surfaces, where the curvature term arises because $\bV$ is not tangent to $\Gamma$.

\subsection{Necessary Optimality Conditions}\label{sec:NoC}
We revisit the necessary optimality conditions of~\eqref{eq:main} from~\cite{HornungMittelmann1990} using shape calculus, in comparison to a differential geometric perspective used in the former. The Lagrangian of~\eqref{eq:main} is given by
\begin{align*}
  &\mathcal{L}(\Omega, \lambda_\text{vol}, \lambda_{i,j})\\
   :=& S(\Gamma_{\text{LA}}) - \left(\sum_{i=1}^\ell \beta_i S(\Gamma_{\text{LS}_i}) \right) + b G(\Gamma) + \lambda_\text{Vol}\left(F(\Gamma) - V_0\right) + \sum_{i=1}^\ell\sum_{j=1}^{k_i} \lambda_{i,j} c(s_{i,j}, \Gamma_{\text{S}_i}).
\end{align*}
We now discuss the shape derivative of the Lagrangian term by term. Because the integrand is a constant with respect to shape perturbations, the shape derivative of most terms in the Lagrangian can immediately be found to be
\begin{equation}\label{eq:AllDerivatives}
\begin{aligned}
d S(\Gamma)[\bV] &= \int\limits_{\Gamma} \div[\Gamma] \bV \dd s = \int\limits_{\Gamma} \inner{\bV}{\bn} \kappa \dd s + \int\limits_{\partial \Gamma} \inner{\bV}{\bmu} \dd \ell \\
d G(\Gamma)[\bV] &= \int\limits_{\Omega} \div[](\inner{\bg}{x}\bV) \dd x \\
&= \int\limits_{\Gamma} \inner{\tilde{\bg}}{\bn}\div[\Gamma]\bV + \inner{[g_is_i\delta_{ij}]_{ij}\bV}{\bn} + \inner{\tilde{\bg}}{d\bn[\bV]} \dd s = \int\limits_{\Gamma}\inner{\bV}{\bn} \inner{\bg}{s} \dd s \\
d F(\Gamma)[\bV] &= \int\limits_{\Omega} \div[] \bV \dd x = \int\limits_{\Gamma} \frac{1}{3}\inner{s}{\bn}\div[\Gamma] \bV + \frac{1}{3}\inner{\bV}{\bn} + \frac{1}{3}\inner{s}{d\bn[\bV]} \dd s \\
&= \int\limits_{\Gamma} \inner{\bV}{\bn} \dd s,
\end{aligned}
\end{equation}
where $\delta_{ij}$ is the Kroenecker delta. The distance function $\dist_{\Gamma_{\text{LS}_i}}(s)$ is independent of the deformation parameter $\varepsilon$. Hence, the shape derivative vanishes and the material derivative is given by transport only, i.e.,
\begin{align*}
d\dist_{\Gamma_{\text{LS}_i}}(s)[\bV] = \dD\dist_{\Gamma_{\text{LS}_i}}(s)\bV(s) =: \inner{\bd_i(s)}{\bV(s)},
\end{align*}
where $\dD$ is the classical spatial Jacobian and $\bd_i(s)$ is the spatial gradient of the distance field with respect to $\Gamma_{\text{LS}_i}$ at point $s$. Picking the respective curvature free expressions from~\eqref{eq:AllDerivatives} with the least regularity requirements on $\Gamma$, one arrives at an alternative variational formulation of~\eqref{eq:YLE}, well-suited for discretization via triangular meshes. Indeed, the following section shows how the respective high-regularity expressions from~\eqref{eq:AllDerivatives} create~\eqref{eq:YLE}.

\subsection{Tangential Motion and the Contact Angle}\label{sec:YLE}
Omitting the subset constraint, one arrives at
\begin{align*}
  &d\mathcal{L}(\Omega, \lambda_\text{vol})[\bV]
  = \int\limits_{\Gamma_{\text{LA}}} \inner{\bV}{\bn} \kappa \dd s + \sum_{i=1}^k \int\limits_{\partial \Gamma_i} \inner{\bV}{\bmu_\text{A} - \beta_i \bmu_i} \dd \ell\\
   &\hspace{5cm}+ b\int\limits_{\Gamma} \inner{\bV}{\bn} \inner{\bg}{s} \dd s + \lambda_{\text{Vol}} \int\limits_{\Gamma} \inner{\bV}{\bn} \dd s,
\end{align*}
where $\partial \Gamma_i$ is the contact line between the liquid to air interface and the contact surface on solid body $i$. In the absence of gravity, i.e. $b=0$, the above expression can only vanish, if the liquid air interface $\Gamma_{\text{LA}}$ is a surface of constant mean curvature, i.e., $\kappa = -\lambda_{\text{Vol}}$, which is indeed the case where~\eqref{eq:YLE} is spatially constant. Furthermore,
\begin{align*}
  \inner{\bV}{\bmu_\text{A} - \beta_i \bmu_i} = 0
\end{align*}
has to hold on $\partial \Gamma_i$, which is not possible for arbitrary perturbation fields $\bV$. If $\Gamma$ is feasible, then $\bmu_i$ is in the tangent space of $\Gamma_{\text{S}_i}$. Hence, if one restricts the motion to tangent directions only, i.e. $\bV = \alpha(s) \bmu_i(s)$, then the necessary optimality condition becomes
\begin{align}\label{eq:TangentMotion1}
 0 = \inner{\bV}{\bmu_\text{A} - \beta_i \bmu_i} = \alpha(s) \inner{\bmu_i}{\bmu_\text{A} - \beta_i \bmu_i} \Leftrightarrow  \inner{\bmu_i}{\bmu_{\text{A}}} - \beta_i = 0.
\end{align}
Because the co-normals need to satisfy $\inner{\bmu_i}{\bmu_{\text{A}}} = \beta_i$ and have unit length, the critical surface must form a contact angle of $\cos \beta_i$ and stationarity is only possible for tangent motions.

\section{Shape Hessians}\label{sec:ShapeHessian}
\subsection{Shape Hessian for Boundary Integrals}
Our desired computational approach is a Newton-iteration to find the roots of $d\mathcal{L}(\Omega, \lambda_\text{vol}, \lambda_{i,j})[\bV, d\lambda_\text{vol}[\bV], d\lambda_{i,j}[\bV]]$. Within the optimization context, this is also called sequential quadratic programming (SQP), because during each update, the second order Taylor-expansion of $\mathcal{L}$ is minimized. For the problem at hand, this means we require the second order shape derivative of functionals of the types
\begin{align*}
	J_1(\Gamma) := \int\limits_{\Gamma} f_1(s)\dd s \text{ and } J_2(\Gamma) := \int\limits_{\Gamma} \inner{\bf_2}{\bn} \dd s,
\end{align*}
the latter arising due to our reformulation of volume integrals to surfaces in~\eqref{eq:S} and~\eqref{eq:G}.

Indeed, it is possible to transform $J_2$ back into a volume integral, which would, however, result in two problems: First, using the volume formulation of the Hessian over $\Omega$ necessitates using a tetrahedral mesh, resulting in considerable numerical overhead. Furthermore, the Hessian matrix post discretization would then suffer from an additional rank deficit for vertices in the interior. Second, using a boundary formulation of the volumetric problem description requires the computation of the total curvature $\kappa$, which is an edge concentrated quantity for the shell meshes consisting of planar triangles our code is using. Hence, optimization and discretization would not commute in this situation.

A repeated application of~\eqref{eq:SurfDerive1} leads to
\begin{equation}\label{eq:ShapeHess1}
\begin{aligned}
&d^2J_1(\Gamma)[\bV, \bW]\\
=& \int\limits_{\Gamma} f_1\div[\Gamma] \bV\div[\Gamma] \bW + df_1[\bV]\div[\Gamma] \bW + d(f_1\div[\Gamma] \bV + df_1[\bV])[\bW]  \ \dd s\\
=& \int\limits_{\Gamma} f_1\div[\Gamma] \bV\div[\Gamma] \bW + df_1[\bV]\div[\Gamma] \bW\\
&\hspace{0.5cm} + df_1[\bW]\div[\Gamma] \bV + f_1d(\div[\Gamma] \bV)[\bW]+ d^2f_1[\bV, \bW]  \ \dd s.
\end{aligned}
\end{equation}
The above expression is more involved than those readily available in the literature, as the consideration of surface integrals here requires the computation of $d(\div[\Gamma] \bV)[\bW]$, the material derivative of the tangent divergence, whereas the volume situation is usually considered otherwise. Material derivatives commute with algebraic operations, such that
\begin{align*}
  d(\div[\Gamma]\bV)[\bW] = d(\tr(\dD_\Gamma \bV))[\bW] = \tr(d(\dD_\Gamma \bV)[\bW]). 
\end{align*}
Due to~\eqref{eq:MaterialSpatial}, one arrives at
\begin{equation}\label{eq:dDV}
\begin{aligned}
  &d(\dD_\Gamma \bV)[\bW]\\
  =& \dD(d\bV[\bW]) - \dD\bV\dD\bW - d(\dD\bV\bn\bn^T)[\bW]\\
  =& \dD_\Gamma(d\bV[\bW]) - \dD\bV\dD\bW + \dD\bV\dD\bW\bn\bn^T - \dD\bV d\bn[\bW]\bn^T - \dD\bV\bn d\bn[\bW]^T\\
  =& \dD_\Gamma(d\bV[\bW]) - \dD\bV\left( \dD_\Gamma \bW - (\dD_\Gamma \bW)^T\bn\bn^T - \bn\bn^T\dD_\Gamma \bW\right)\\
  =& \dD_\Gamma(d\bV[\bW]) - \dD_\Gamma\bV\dD_\Gamma\bW + \dD\bV\left(\dD_\Gamma \bW\right)^T \bn\bn^T\\
  =& \dD_\Gamma(d\bV[\bW]) - \dD_\Gamma\bV\dD_\Gamma\bW + \dD_\Gamma\bV\left(\dD_\Gamma \bW\right)^T \bn\bn^T.
\end{aligned}
\end{equation}
In the last step, the remaining full Jacobian $\dD\bV$ can be replaced with the tangential Jacobian, because
\begin{align*}
  \dD\bV(\dD_\Gamma\bW)^T\bn\bn^T = (\dD_\Gamma\bV + \dD\bV\bn\bn^T)(\dD_\Gamma\bW)^T\bn\bn^T = \dD_\Gamma\bV(\dD_\Gamma\bW)^T\bn\bn^T.
\end{align*}
Hence, taking the trace results in
\begin{align}\label{eq:MatDiv}
  d(\div[\Gamma]\bV)[\bW] = -\tr\left(\dD_\Gamma\bV\dD_\Gamma\bW\right) + \inner{\left(\dD_\Gamma\bV\right)^T\bn}{\left(\dD_\Gamma\bW\right)^T\bn}.
\end{align}
To achieve symmetry, we assume $d\bV[\bW]=0$ here and in the following. This makes the shape Hessian derived via repeated differentiation align with the expressions obtained by enforcing a vector space structure. It is worth mentioning that the latter part of~\eqref{eq:MatDiv} does not appear when volume integrals are considered.
Summarizing the above, the shape Hessian of $S(\Gamma)$ from~\eqref{eq:S} is given by
\begin{align*}
&d^2 S(\Gamma)[\bV, \bW]\\
 =& \int\limits_\Gamma \div[\Gamma] \bV\div[\Gamma] \bW -\tr\left(\dD_\Gamma\bV\dD_\Gamma\bW\right) + \inner{\left(\dD_\Gamma\bV\right)^T\bn}{\left(\dD_\Gamma\bW\right)^T\bn} \ \dd s.
\end{align*}
\subsection{Shape Hessian of the Normal}
The shape Hessian of $F(\Gamma)$ and $G(\Gamma)$ in~\eqref{eq:S} and~\eqref{eq:G} is more involved, as the second material derivative $d^2\bn[\bV,\bW]$ of the normal is needed. For the first Eulerian derivative, Equation~\eqref{eq:SurfDerive1} is applied to the $J_2$ case above, leading to
\begin{align}\label{eq:ShapeGrad2}
  dJ_2(\Gamma)[\bV] = \int\limits_{\Gamma} \inner{\bf_2}{\bn}\div[\Gamma]\bV + \inner{d\bf_2[\bV]}{\bn} + \inner{\bf_2}{d\bn[\bV]} \ \dd s
\end{align}
and the symmetric repeated differentiation second order shape derivative is given by
\begin{equation}\label{eq:ShapeHessNormal}
\begin{aligned}
  &d^2J_2(\Gamma)[\bV, \bW]\\
  =& \int\limits_{\Gamma} \inner{\bf_2}{\bn}\div[\Gamma]\bV\div[\Gamma]\bW -\inner{\bf_2}{\bn}\tr\left(\dD_\Gamma\bV\dD_\Gamma\bW\right) \\
  &\hspace{0.2cm} +\inner{\bf_2}{\bn}\inner{\left(\dD_\Gamma\bV\right)^T\bn}{\left(\dD_\Gamma\bW\right)^T\bn} \\
  &\hspace{0.2cm}+\left(\inner{d\bf_2[\bV]}{\bn} + \inner{\bf_2}{d\bn[\bV]}\right)\div[\Gamma]\bW+\left(\inner{d\bf_2[\bW]}{\bn} + \inner{\bf_2}{d\bn[\bW]}\right)\div[\Gamma]\bV\\
  &\hspace{0.2cm}+\inner{d\bf_2[\bV]}{d\bn[\bW]}+\inner{d\bf_2[\bW]}{d\bn[\bV]}\\
  &\hspace{0.2cm}+ \inner{d^2\bf_2[\bV,\bW]}{\bn} + \inner{\bf_2}{d^2\bn[\bV, \bW]} \ \dd s.
\end{aligned}
\end{equation}
The second material derivative of the normal can be explicitly computed using~\eqref{eq:dn},
\begin{equation}\label{eq:ddn}
\begin{aligned}
&d^2\bn[\bV, \bW]\\
 =& d(d\bn[\bV])[\bW] = -d(\dD_\Gamma\bV)[\bW]^T\bn - (\dD_\Gamma\bV)^Td\bn[\bW]\\
 \stackrel{\eqref{eq:dDV}}{=}& -(\dD_\Gamma(d\bV[\bW]))^T\bn + (\dD_\Gamma\bW)^T(\dD_\Gamma \bV)^T\bn\\
 &\hspace{3.0cm} - \bn\bn^T(\dD_\Gamma\bW)(\dD_\Gamma\bV)^T\bn + (\dD_\Gamma\bV)^T(\dD_\Gamma\bW)^T\bn\\
 =& (\dD_\Gamma \bW)^T(\dD_\Gamma\bV)^T\bn - \inner{(\dD_\Gamma\bV)^T\bn}{(\dD_\Gamma\bW)^T\bn}\bn + (\dD_\Gamma\bV)^T(\dD_\Gamma\bW)^T\bn,
\end{aligned}
\end{equation}
where we have again used the independency $d\bV[\bW]=0$ in the last step. Thus, when inserting~\eqref{eq:ddn} into~\eqref{eq:ShapeHessNormal}, one arrives at an explicit expression for the second order shape derivative $d^2J_2[\bV,\bW]$ in a general setting, which can now easily be adapted to the second derivative of $F(\Gamma)$ and $G(\Gamma)$ in~\eqref{eq:S} and~\eqref{eq:G}, using the respective material derivatives and the independency $d\bV[\bW]=0$, namely
\begin{align*}
  \bf_2(s) = \frac{1}{3}s, \quad d\bf_2[\bV](s) = \frac{1}{3}\bV(s), \quad d^2\bf_2[\bV, \bW](s) = \frac{1}{3}d\bV[\bW](s) \equiv 0
\end{align*}
for $F(\Gamma)$ and
\begin{align*}
  \bf_2(s) = \frac{1}{2}\left(g_is_i^2\right)_{i=1}^3, \quad d\bf_2[\bV](s) = \left(g_is_iV_i\right)_{i=1}^3, \quad d^2\bf_2[\bV, \bW](s) = \left(g_iV_iW_i\right)_{i=1}^3
\end{align*}
for $G(\Gamma)$. Finally, the second derivative of the $\Gamma$-independent distance-function is readily given by
\begin{align*}
  d^2\dist_{\Gamma_{\text{LS}_i}}[\bV, \bW] = \bV^T\bH\bW,
\end{align*}
where $\bH$ is the spatial Hessian matrix of the distance field.

\section{Numerical Implementation}\label{sec:Implementation}
\subsection{Shape Hessian and Finite Elements}
To numerically find the shape of the capillary bridge, we use finite elements on a shell mesh $\Gamma_h$ of topological dimension two consisting of planar triangles $T$ in $\R^3$ provided by the FEniCS framework~\cite{Fenics_book}. The classical continuous Galerkin finite element space of order $r$ of $d$-dimensional vectors is given by
\begin{align*}
  \mathcal{CG}_r^d(\Gamma) := \left\{\bu\in \left(L^2(\Gamma_h)\right)^d : \bu_{\vert T} \in \mathcal{P}^d_r(T)\right\},
\end{align*}
where $\mathcal{P}_r$ is the space of polynomials of order $r$. To include the respective scalar constraints of volume and distance, we augment this space to $Q :=  \mathcal{CG}_1^3(\Gamma) \times\R\times\R^{n_in_j}$ to include the adjoint multipliers. The Newton update for iteratively finding the capillary bridge is then interpretable as a variational problem, namely to find $q[\bW] := (\bW, d\lambda_\text{vol}[\bW], d\lambda_{i,j}[\bW]) \in Q$, such that
\begin{align}\label{eq:KKT}
  d^2\mathcal{L}(\Gamma_h^{k}, \lambda_\text{vol}^{k}, \lambda_{i,j}^k)[\bV, \bW] = -d\mathcal{L}(\Gamma_h^k, {\lambda}_\text{vol}^k, {\lambda}_{i,j}^k)[\bV] \quad \forall q[\bV] \in Q.
\end{align}
Post discretization, this variational problem creates the typical Karush-Kuhn-Tucker (KKT) system. Finally, the state is updated via
\begin{equation*}
(\Gamma^{k+1}_h, \lambda_\text{vol}^{k+1}, \lambda_{i,j}^{k+1}) := (\Gamma^{k}_h, \lambda_\text{vol}^{k}, \lambda_{i,j}^{k}) + q[\bW].
\end{equation*}
We use the discrete shape differentiability capabilities of FEniCS~\cite{FenicsShapeAD} to numerically validate the $\bV$-$\bW$ block in~\eqref{eq:KKT}, confirming that indeed for our first and second order shape derivatives the ``discretize-then-optimize'' and ``optimize-then-discretize'' approach commutes irrespective of mesh size. Including the off-diagonal blocks in this confirmation is unfortunately not conveniently possible, as FEniCS does not provide a means to compute the $\bV$-$d\lambda[\bW]$ off-diagonals discretely.

A straightforward implementation of~\eqref{eq:KKT} results in the challenge of globalizing the convergence, i.e., guaranteeing a positive definite Hessian matrix away from the optimum. To address this, we have also implemented an approximative Newton scheme, where the $\bV$, $\bW$ block in~\eqref{eq:KKT} is replaced by the $H^1(\Gamma)$ scalar product
\begin{align}\label{eq:H1Smooth}
  \inner{\bV}{\bW}_{H^1} := \int\limits_{\Gamma} \gamma\inner{\nabla_\Gamma \bV}{\nabla_\Gamma \bW} + \inner{\bV}{\bW}\ \dd s,
\end{align}
where $\gamma>0$ is some smoothing parameter. This approach is sometimes also called a Sobolev gradient descent and the resulting KKT-system is positive definite.

There are several reasons for rank deficiencies of the proper KKT system for shapes, also at the optimum. Using a full 3D deformation $\bV, \bW \in\mathcal{CG}_1^3(\Gamma_h)$ is indeed too rich. As can be seen from~\eqref{eq:SurfDerive2} and~\eqref{eq:SurfDerive3}, tangential motions are in the kernel of~\eqref{eq:KKT}. Furthermore, the second order shape derivatives~\eqref{eq:ShapeHess1}, \eqref{eq:ShapeGrad2} and \eqref{eq:ShapeHessNormal} only involve spatial derivatives of $\bV$ and $\bW$, such that translations are also within the kernel. As such, one typically restricts the deformation unknown to $\bV = \alpha\bN$, where $\alpha \in \mathcal{CG}_1^1(\Gamma_h)$ and $\bN$ is some vertex average of the cell normal $\bn$. However, as discussed in Section~\ref{sec:NoC}, the contact angle requirement~\eqref{eq:TangentMotion1} demands a tangential motion of the vertices constituting $\Gamma_{\text{LS}_i}$ and in particular $\partial \Gamma_i$. As such, we restrict the motion fields
\begin{align}\label{eq:NewtonRank1}
\bV(s_i) =
\begin{cases}
  \alpha_1(s_i)\bmu_\text{LA}(s_i) + \alpha_2(s_i)\bmu_\text{LS}(s_i), \quad s_i \in \Gamma_\text{LA} \cap \Gamma_{\text{LS}}\\
  \alpha(s_i)_1\bN(s_i), \quad \text{otherwise},
\end{cases}
\end{align}
where $\alpha_1$, $\alpha_2$ are $\mathcal{CG}_1^1(\Gamma_h)$ functions, likewise for $\bW$. Thus, we admit a motion in normal direction only for the vertices, except for those forming the solid-gas-liquid triple phase lines. There, each vertex is allowed to move in their respective two dimensional subspace spanned locally by the two edge-averaged co-normals $\bmu_\text{LA}$ and $\bmu_\text{LS}$. It is worth mentioning that as a consequence, once feasibility with respect to the distance constraint is achieved, there is no longer any possible motion of the interior vertices of $\Gamma_\text{LS}$, except for the contact line $\partial \Gamma := \Gamma_\text{LA} \cap \Gamma_{\text{LS}}$, which can still move in the plane spanned locally by the vertex co-normals. Hence, if one switches to the Newton-solver too early, the mesh quality and point spacing of $\Gamma_\text{LS}$ is going to degrade. At present, we counter this problem by remeshing, which is necessary as we aim at being able to deal with large scale deformations from the initial guess, anyway. Incorporating equal vertex spacing into the optimization goal to remove the tangential kernel of the Hessian for vertices in the interior of $\Gamma_\text{LS}$ is part of future work.

Depending on the geometrical setup of the solid bodies $\Gamma_{\text{S}_i}$, translations can also be in the kernel of~\eqref{eq:KKT}. One such case would be a capillary bridge between two planes. Full rank in this situation is achieved by activating additional centroid constraints. Surface and volume centroid are given by
\begin{align*}
  \partial c_i(\Gamma) := \frac{1}{S(\Gamma)} \int\limits_\Gamma s_i \dd s,\quad c_i(\Omega) := \frac{1}{\text{vol}(\Omega)} \int\limits_\Omega x_i \dd x.
\end{align*}
While using the surface centroid might seem to be the approach of choice for a code operating on the shell alone, we rather choose the volume centroid, as the expression can be considerably simplified by exploiting the volume constraint in this setting. Thus, we constrain the $i$-th component of the centroid via
\begin{align*}
 \tilde{x}_i \stackrel{!}{=} \frac{1}{\text{vol}(\Omega^*)}\int\limits_\Omega x_i \dd x = \frac{1}{\text{vol}(\Omega^*)}\int\limits_\Gamma \inner{\bf_2^i}{\bn}\dd s,
\end{align*}
where $\text{vol}(\Omega^*)$ is the fixed target volume, $\bf_2^i(s) := \frac{1}{2}\left(s_{j}^2\delta_{ij}\right)_{j=1}^3$ and $\delta_{ij}$ is the Kronecker-Delta. The respective first and second order Eulerian derivatives are again given by~\eqref{eq:ShapeGrad2}, \eqref{eq:ShapeHessNormal} and~\eqref{eq:ddn}, where
\begin{align*}
d\bf_2^i[\bV] = (V_j\delta_{ij})_{j=1}^3, \quad d^2\bf_2^i[\bV, \bW] \equiv 0.
\end{align*}
More details on these constraints can also be found in~\cite{SAD_2018}.


\section{Numerical Results}\label{sec:NumericalResults}
\subsection{Sphere-Plate, the unique situation}
In order to validate our second order optimization scheme, we consider the case of a spherical body with zero distance from a plane, as this situation is very well understood theoretically. In the following, we base our evaluation on the data available in~\cite{orr_scriven_rivas_1975}. In particular, there is a closed form solution of the total force available as
\begin{align}\label{eq:ForceOrr}
  F_t = 2\pi\sigma R\left[\sin \psi \sin(\theta_1 + \psi) - HR\sin^2 \psi\right],
\end{align}
where $\theta_i$ are the wetting angles at the sphere and plane. Furthermore, $\psi$ is the so called filling angle. The available data is given in non-dimensionalized form as the ratio of contact force $F_t$ to radius $R$ times surface tension $\sigma$. However, there is no closed form solution for the curvature $H$ as a coupled system of non-linear equations has to be solved. Tabularized values are provided in~\cite{orr_scriven_rivas_1975} depending on the angles $\theta_1$ and $\psi$ for the curvature. Because our code requires the desired contact angle and volume as input rather than the angle $\theta_1$, we use a volume constraint of $0.165$ as per the tabularized values in~\cite{orr_scriven_rivas_1975}, corresponding to the angles $\theta_1 = \theta_2 = 40^\circ$. Seeing that curvature is not readily available for a discrete surface of planar triangles, we use the adjoint of the volume constraint for comparisons, utilizing the correspondence
\begin{align*}
  \Delta p = 2HR = -\lambda_\text{vol}
\end{align*}
as described in Section~\ref{sec:YLE} and validate the resulting effective pressure difference $\Delta p$ as the corresponding quantity. Instead of using the close-form representation of the forces~\eqref{eq:ForceOrr} for spheres, we follow the general expression for forces~\cite{virozub_2009} and compute the force acting on obstacle $\Gamma_{\text{S}_i}$ via numerical quadrature of the following integrals
\begin{align}\label{eq:Forces}
  \bF_{p,i} := \lambda_\text{vol}\int\limits_{\Gamma_{\text{S}_i}} \bn \ \dd s \in \R^3, \quad \bF_{s,i} := \int\limits_{\partial\Gamma_i} \bmu\ \dd \ell \in \R^3,
\end{align}
with $\bF_{t,i} := \bF_{p,i} + \bF_{s,i}$ being the resulting force acting on obstacle $i$. The vector $\bmu$ is again the co-normal on $\partial\Gamma_i$, e.g. a vector orthogonal to $\partial\Gamma_i$ and tangential to $\Gamma_\text{LA}$. In particular the integrands in $\bF_p$ and $\bF_s$ are constants per cell or edge for our mesh consisting of planar triangles. As such, the above integrals can readily be computed by summing the contribution of each triangle and edge individually.

The shapes of the two obstacles, the sphere and plane, enters our code via closed form descriptions for the signed distances in~\eqref{eq:DistanceConstraints}, in particular
\begin{align*}
  \dist(x, \Gamma_{\text{S}_1}) &= \sqrt{x_1^2 + x_2^2 + (x_3-z)^2} - R,\\
  \dist(x, \Gamma_{\text{S}_2}) &= x_3 + z,
\end{align*}
where $z=0$ is twice the distance between the objects, i.e., $z=0$ here. First and second order spatial derivatives of these functions, which are needed for the shape derivative and shape Hessian, can easily be computed and the non-differentiability of the square root at zero is irrelevant, as the centroid of the sphere is outside of our computational mesh $\Gamma_h$.

\begin{figure}[h!]
\begin{center}
\includegraphics[width=0.45\textwidth]{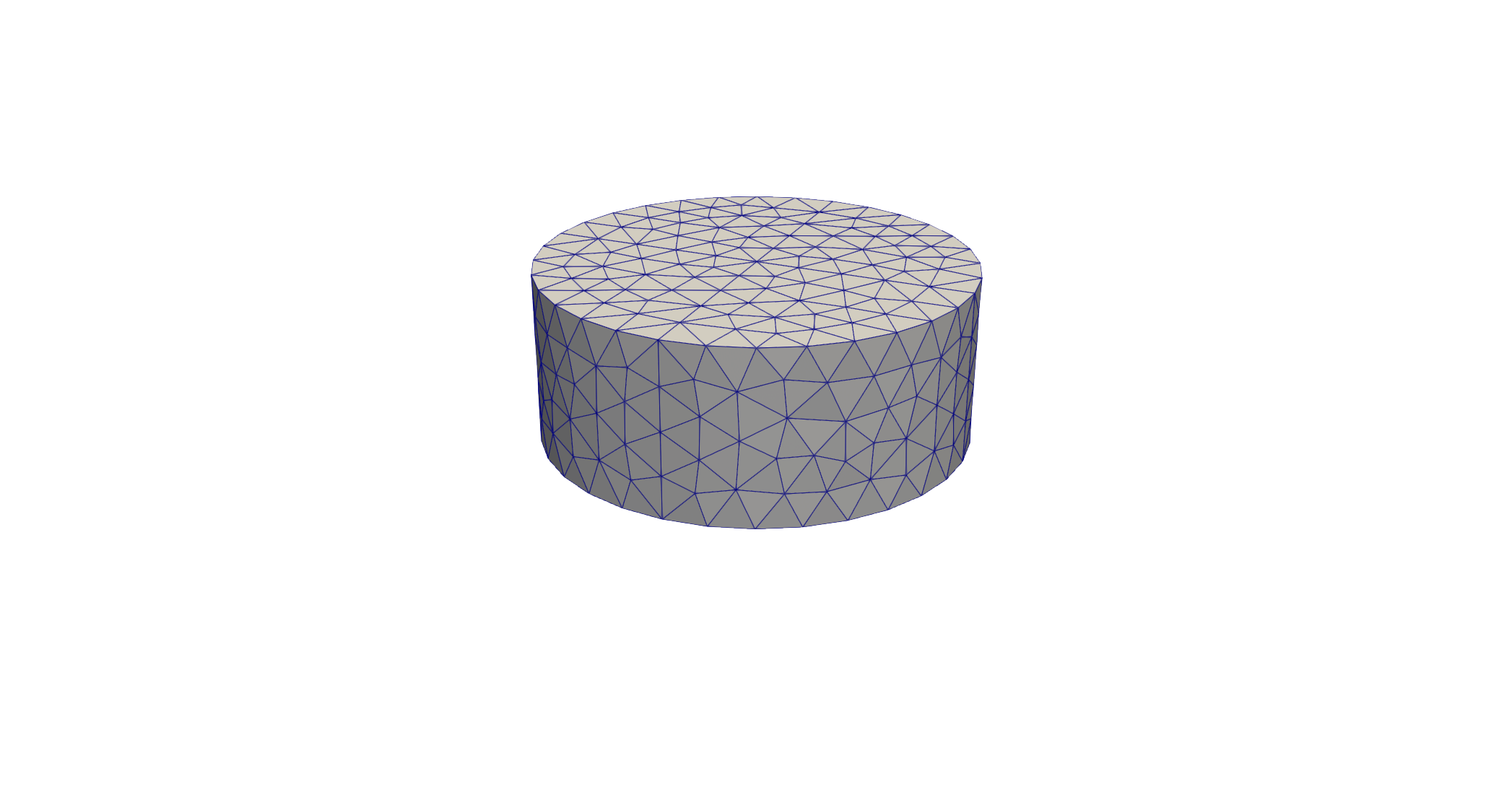}
\includegraphics[width=0.45\textwidth]{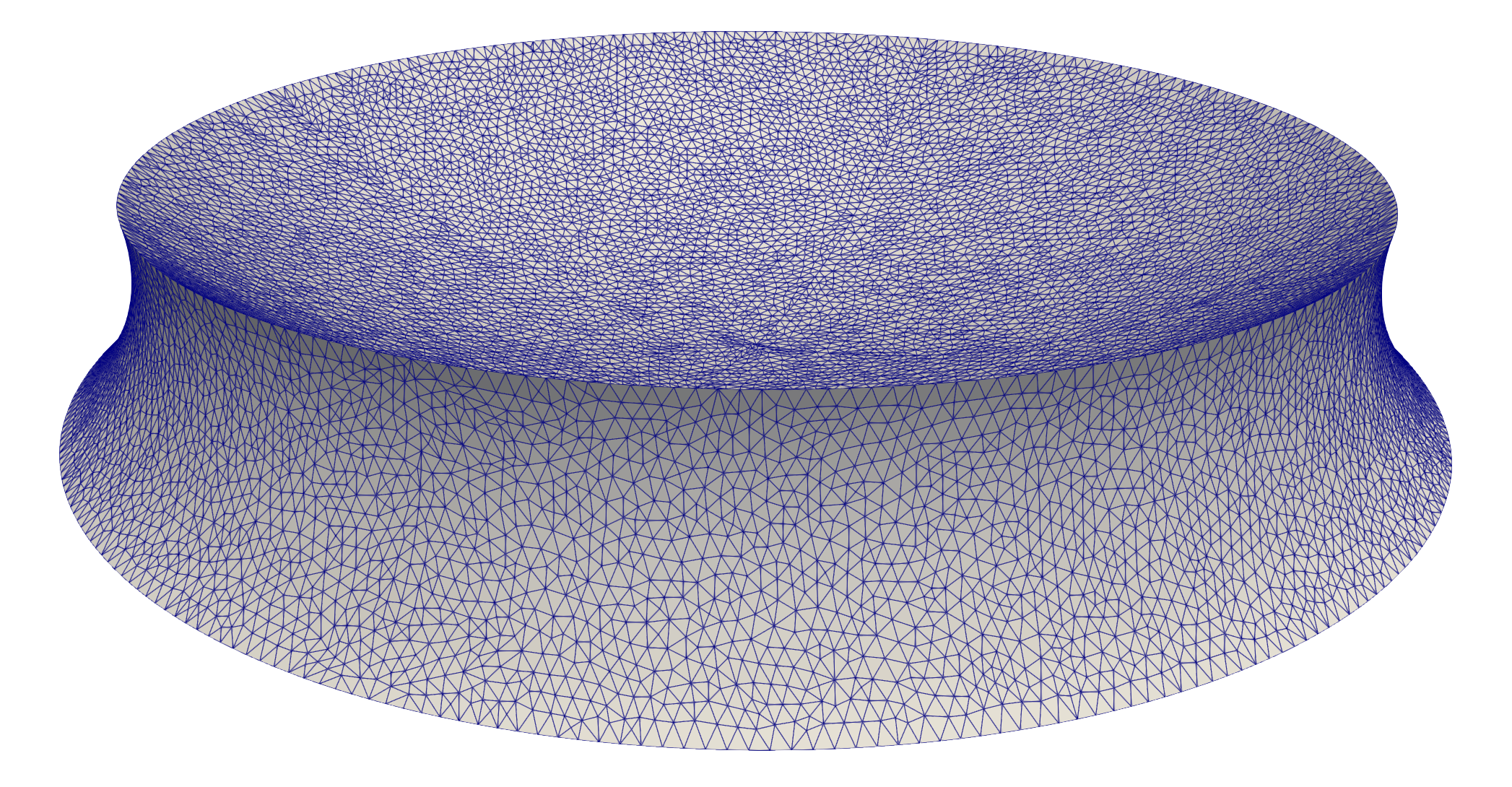}
\end{center}
\caption{Initial cylinder ($3,686$ triangles) and final shape ($52,440$ triangles) for the comparison with~\cite{orr_scriven_rivas_1975}.}\label{fig:OrrResult}
\end{figure}
Our initial geometry is a cylinder around the $x_3$-axis of radius $1.0$, spanning from $x_3 = -0.1$ to $x_3 = 0.1$, consisting of $3,686$ triangles. Hence, we start infeasible with respect to all constraints. The top cap of the cylinders is subject to the touch condition based on $\dist(x, \Gamma_{\text{S}_1})$, while the bottom cap has to adhere $\dist(x, \Gamma_{\text{S}_2})=0$ as a constraint. We first conduct $5$ steps using the approximate Newton-Scheme, where the second order partial shape derivative in the KKT-System, i.e., the corresponding diagonal block, is approximated with the $H^1$-inner product~\eqref{eq:H1Smooth} with smoothing $\gamma = 1$.
After $5$ steps, we remesh $\Gamma_h$ using the compounding of gmsh 3.0.6~\cite{gmsh_remesh2, gmsh_remesh1}. After the remeshing, we switch to the proper Newton iteration. As the geometric situation under consideration is rotationally invariant, we also have to restrict the motion of the vertices once we move to the Newton scheme as discussed in~\eqref{eq:NewtonRank1}. As a consequence, interior vertices of the cap and bottom of the original cylinder, or their descendants after remeshing, can no longer slide over the obstacles once we move to the Newton scheme. However, we found that after the first initial $5$ gradient steps, where sliding is possible, the unknown surface $\Gamma$ already clings to the solid bodies in a well-behaved manner.
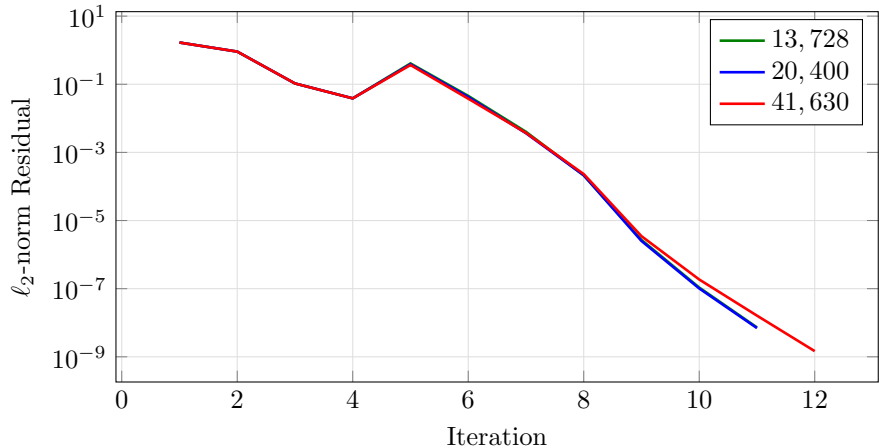
\begin{figure}[h!]
\begin{center}
\begin{tikzpicture}
\begin{axis}[
	ymode=log,
	xlabel=Iteration, 
	ylabel=$\ell_2$-norm Residual,
	grid=both,
	minor grid style={gray!25},
	major grid style={gray!25},
	width=0.9\linewidth,
	height=0.5\linewidth,
	no marks]
\addplot[line width=1pt,solid,color=green!50!black] table[x index = {0}, y index = {2},col sep=space]{Data/meshresolution/orrStandard_0.03/history.txt};\addlegendentry{$13,728$};
\addplot[line width=1pt,solid,color=blue] table[x index = {0}, y index = {2},col sep=space]{Data/meshresolution/orrStandard_0.024/history.txt};\addlegendentry{$20,400$};
\addplot[line width=1pt,solid,color=red] table[x index = {0}, y index = {2},col sep=space]{Data/meshresolution/orrStandard_0.017/history.txt};\addlegendentry{$41,630$}; 
\end{axis}
\end{tikzpicture}
\end{center}
\caption{Convergence history for the validation with data from~\cite{orr_scriven_rivas_1975} for different number of cells. The Newton iteration is activated after remeshing in iteration 5.}\label{fig:OrrConvergence}
\end{figure}
\begin{table}[h!]
\begin{center}
\begin{tabular}{cllll}
\#Cells & $\Delta p = -\lambda_\text{vol}$ & $\Vert \bF_t\Vert$ & Surface Area & Error in Force\\
\hline
- & $-2.6435$ & $7.4088$ & $1.0990$ & - \\
\hline
$13,728$ & $-2.641839$ & $7.434081$ & $1.096037$ & $+0.3412\%$\\ 
$20,400$ & $-2.641954$ & $7.432478$ & $1.096901$ & $+0.3196\%$\\ 
$41,630$ & $-2.642561$ & $7.428643$ & $1.097786$ & $+0.2678\%$\\ 
\end{tabular}
\end{center}
\caption{Reference data in comparison to our computations for varying number of triangles. The original tabularized data from~\cite{orr_scriven_rivas_1975} (first line) has been non-dimensionalized as $\Vert \bF_t\Vert/(R\sigma)$ to match our code. Forces on sphere shown.}\label{tab:Orr}
\end{table}
Initial and final mesh are shown in Figure~\ref{fig:OrrResult} and the respective quantities of interest are listed in Tab.~\ref{tab:Orr}. The convergence of the optimization residual is shown in Figure~\ref{fig:OrrConvergence}.

Our methodology achieves very accurate results in comparison to the reference data in~\cite{orr_scriven_rivas_1975}. In particular, the relative error in the computed total force is less than $0.5\%$ for all meshes considered, in particular including the coarsest grid of only $13,728$ cells. We also observe excellent convergence speed of the Newton method, which is activated after one remeshing in iteration 5 for all grids. Observing the convergence history in Figure~\ref{fig:OrrConvergence}, one can see that the reduction of the residual with the Newton method is indeed very rapid, but not fully grid independent. A possible source for this slight mesh dependency can possibly found in the subset constraint, which is implemented discretely on a per-vertex level without considering a possible analytic or function space equivalent of the subset constraint.

\subsection{Sphere-Plate, Non-Axisymmetric Situation}
\begin{figure}[h!]
\begin{center}
\includegraphics[width=0.2\textwidth]{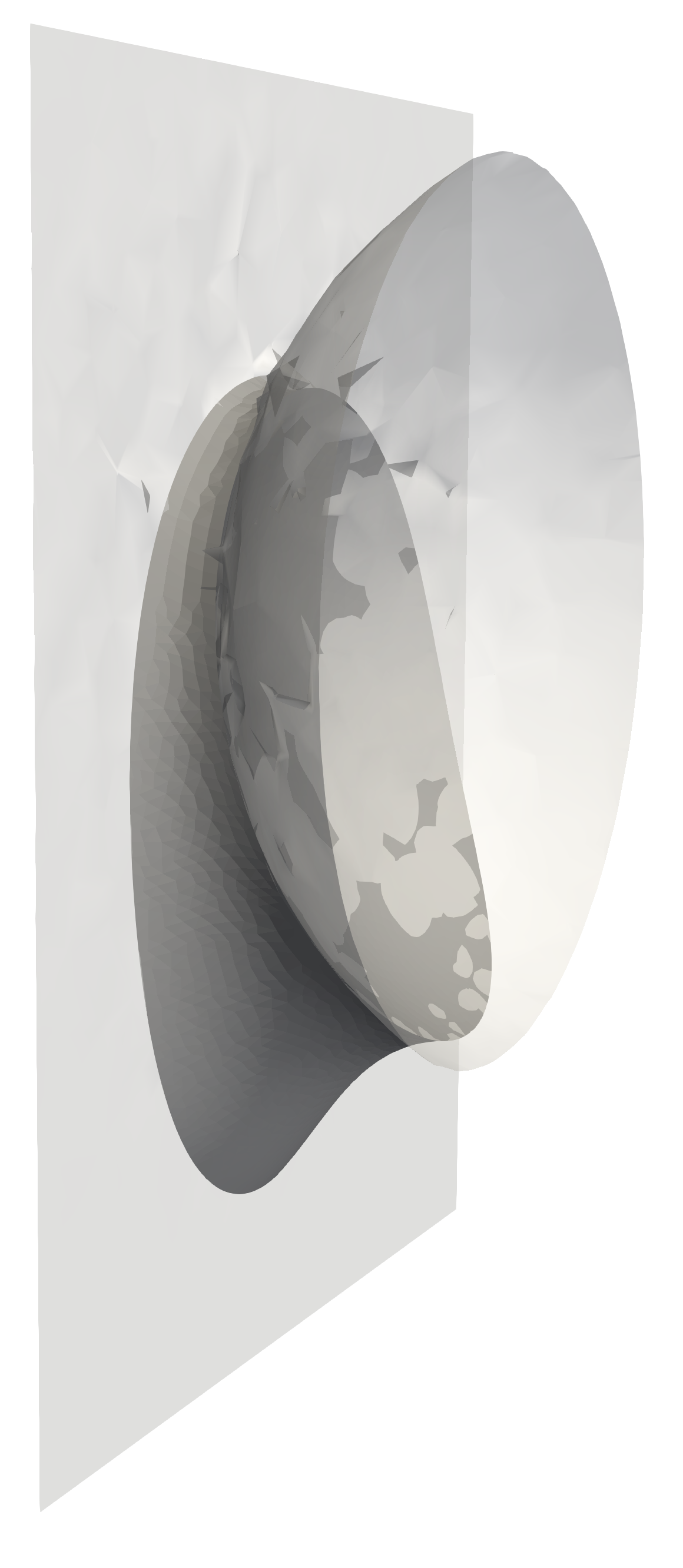}
\includegraphics[width=0.2\textwidth]{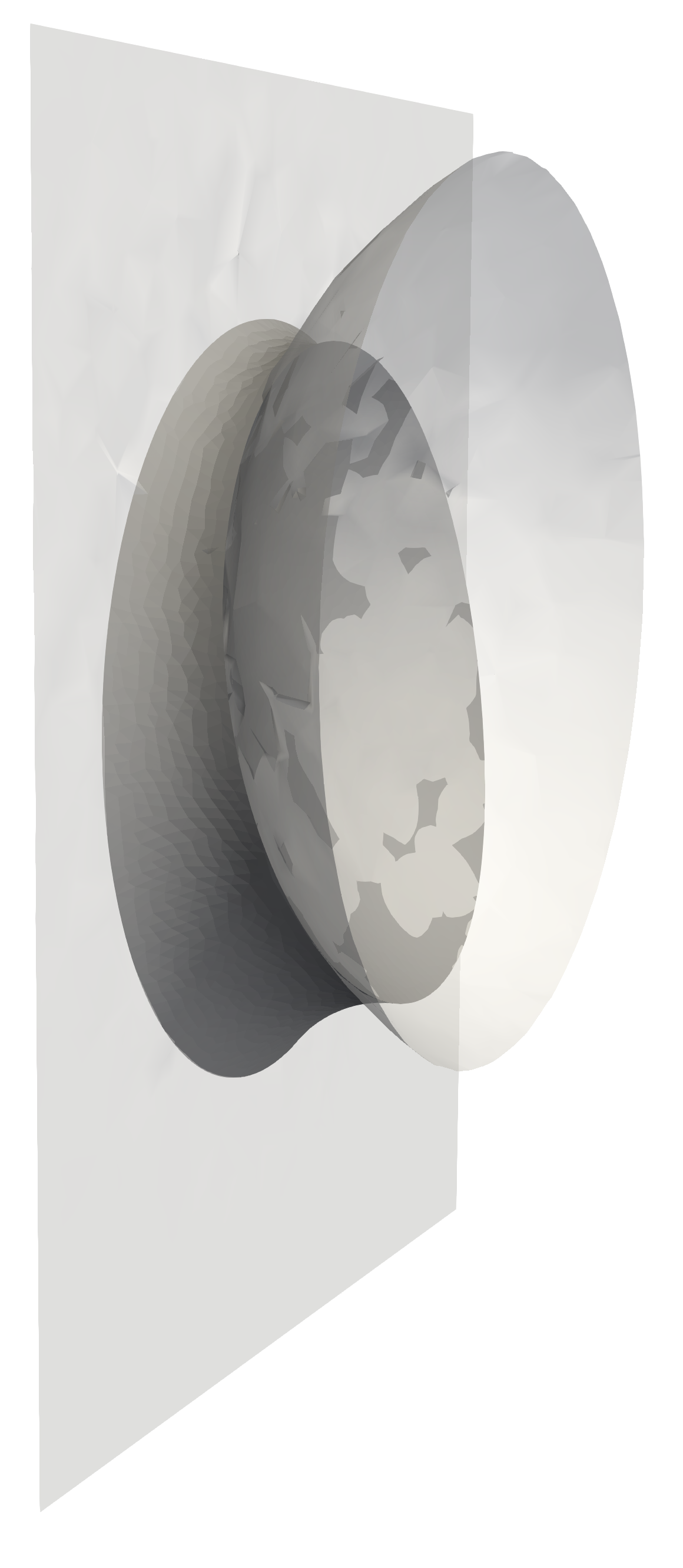}
\includegraphics[width=0.2\textwidth]{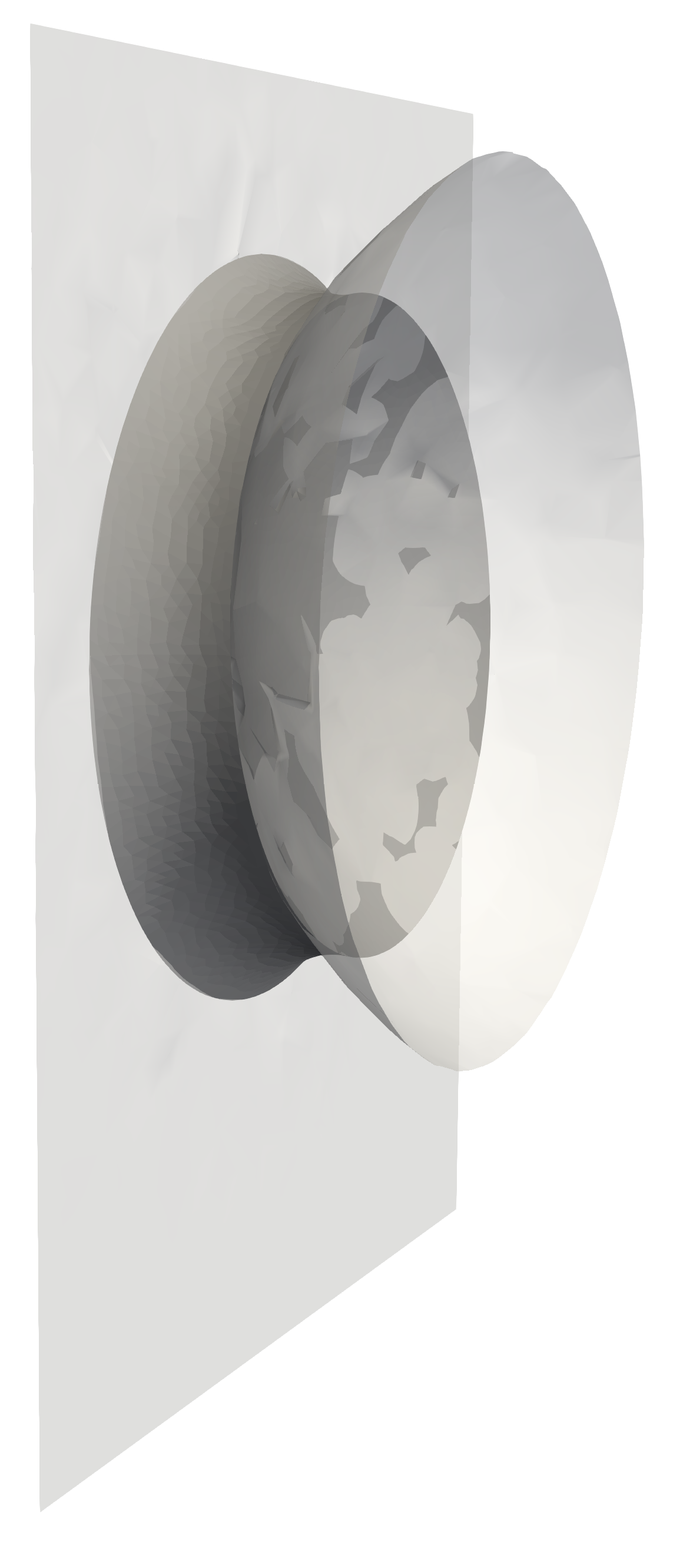}
\includegraphics[width=0.2\textwidth]{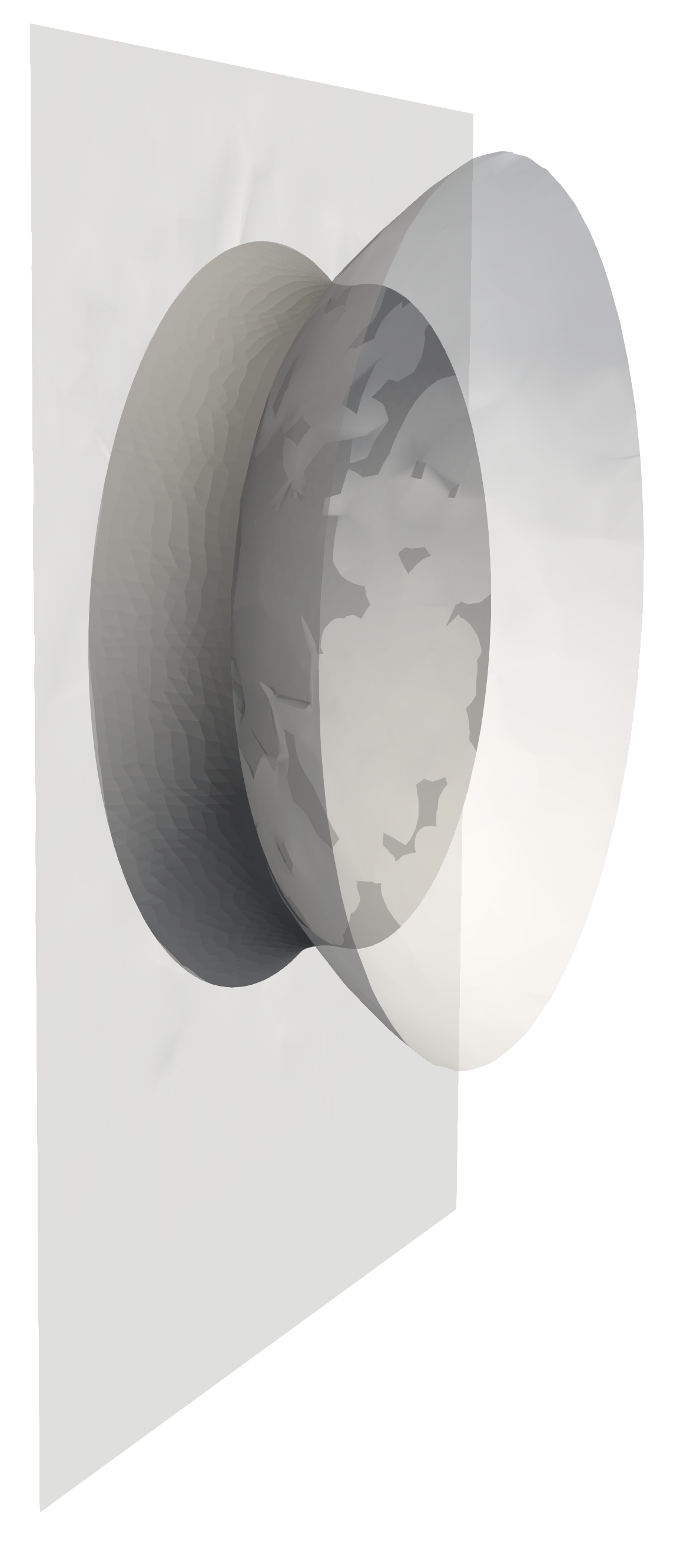}
\caption{The same test case with Bond numbers $8.0$, $4.0$, $1.0$ and $0.5$. Gravity applied in positive $x_2$-direction plotted downwards.}\label{fig:OrrWithGravity}
\end{center}
\end{figure}

We also use the same sphere-plate setting to study the non-axisymmetric situation. To this end, we apply gravity in the negative $x_2$-direction with varying bond numbers of $8.0$, $4.0$, $1.0$ and $0.5$. We did not obtain solutions for bond numbers $10.0$ and up, as these would lead to topological changes with the droplet splitting and detaching in part from the obstacle. The respective shapes are shown in Figure~\ref{fig:OrrWithGravity}.
In particular,~\cite{orr_scriven_rivas_1975} states that the gravitational effect is negligible provided that the Bond number satisfies $b\ll1$. For comparison, we have tabulated the computed force for different Bond numbers in Table~\ref{tab:BondOrr}.
\begin{table}[h!]
\begin{center}
\begin{tabular}{crrrrr}
Bond & $\Delta p=-\lambda_\text{vol}$ & $F_{t,x_1}$  & $F_{t,x_2}$ & $F_{t,x_3}$ & $\Vert \bF_{t}\Vert$\\
\hline
$0.5$ & $-2.6658$ & $ 8.283069\cdot 10^{-5}$ & $-0.150065$ & $7.456180$ & $7.4577$\\
$1.0$ & $-2.7294$ & $ 8.731053\cdot 10^{-4}$ & $-0.301205$ & $7.519072$ & $7.5251$\\
$4.0$ & $-3.9041$ & $-6.616371\cdot 10^{-4}$ & $-1.301884$ & $8.597352$ & $8.6954$\\
$8.0$ & $-7.0027$ & $4.353426\cdot10^{-3}$ & $-2.727429$ & $10.868704$ & $11.2057$
\end{tabular}
\end{center}
\caption{Resulting total forces on the sphere in components and magnitude for different Bond numbers with gravity applied in positive $x_2$-direction for the sphere-plate case from~\cite{orr_scriven_rivas_1975}.}\label{tab:BondOrr}
\end{table}

\subsection{Sphere-Plate, the Non-Unique Situation}
We revisit the situation of a sphere of unit radius $1.0$ placed over a plane with a non-dimensional distance of $2.2$. The desired contact angles are $30^\circ$ on the sphere and $80^\circ$ on the plane. This setting was studied analytically in~\cite{rubinstein_2014}, where four unduloid bridges are given, all satisfying the same non-dimensional volume of $V=3.6$, a constant curvature and the desired touch and angle constraints with the obstacles. Hence, they are critical shapes of~\eqref{eq:main}. However, as mentioned in~\cite{rubinstein_2014}, it is unclear which of these candidates constitutes a physical solution. The respective unduloid bridges do not admit a closed-form description of their geometry. However, the authors of~\cite{rubinstein_2014} kindly provided us with a polygonal graph of two of their unduloids as a cut through the $z=0$ plane: The unduloid with a meniscus of $32^\circ$ was given as a plane graph of $363$ points and the unduloid with a meniscus angle of $47.2^\circ$ was kindly provided as a plane graph of $443$ points. We then use the compounding functionality of gmsh 3.0.6~\cite{gmsh_remesh2, gmsh_remesh1} to create a shell mesh of the resulting bodies of revolution with optimal point spacing independent of the original points provided in the graphs. These two shells are then used as the starting geometry in our program.
\begin{figure}[h!]
\begin{center}
\includegraphics[width=0.32\textwidth]{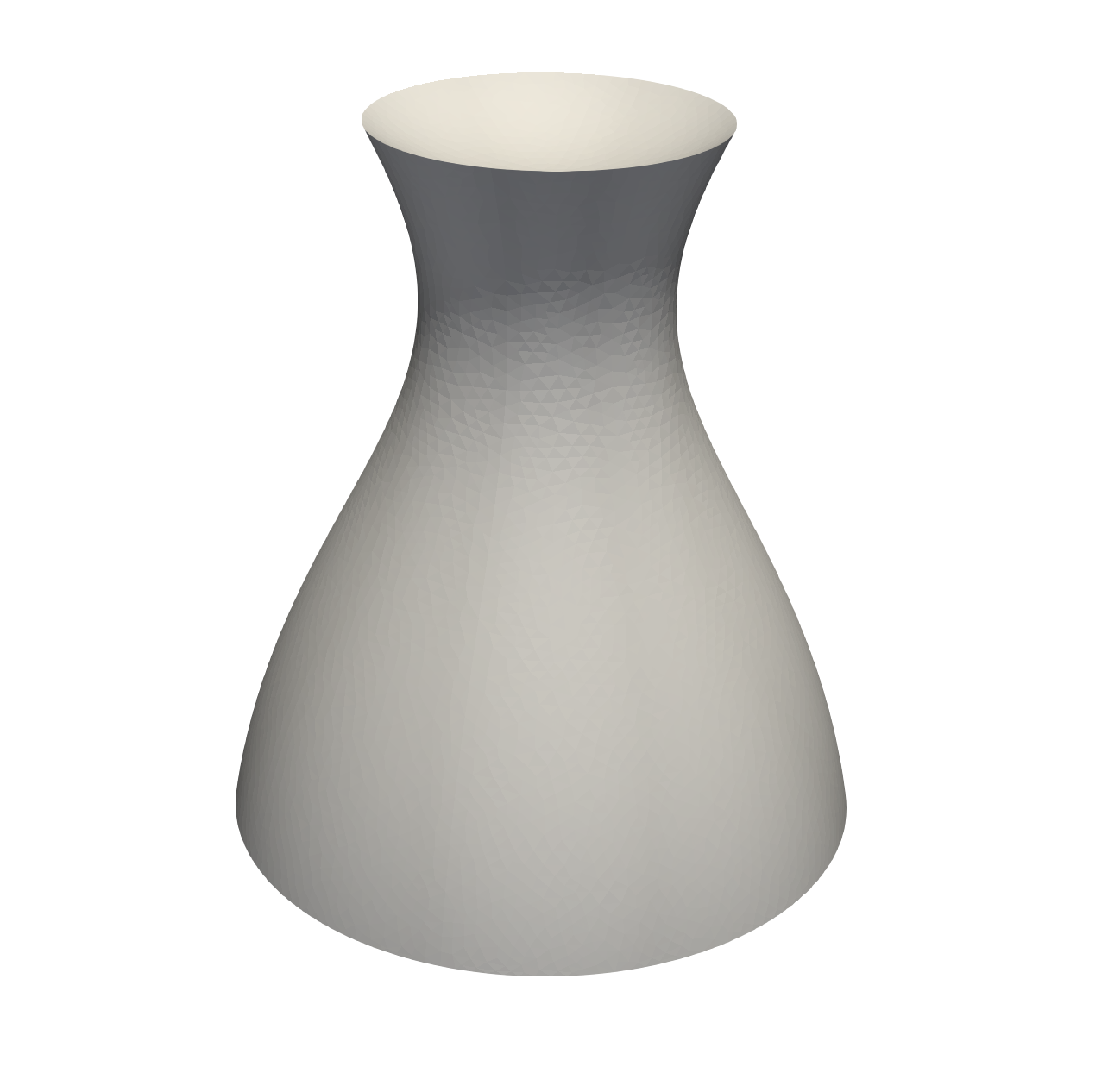}
\includegraphics[width=0.32\textwidth]{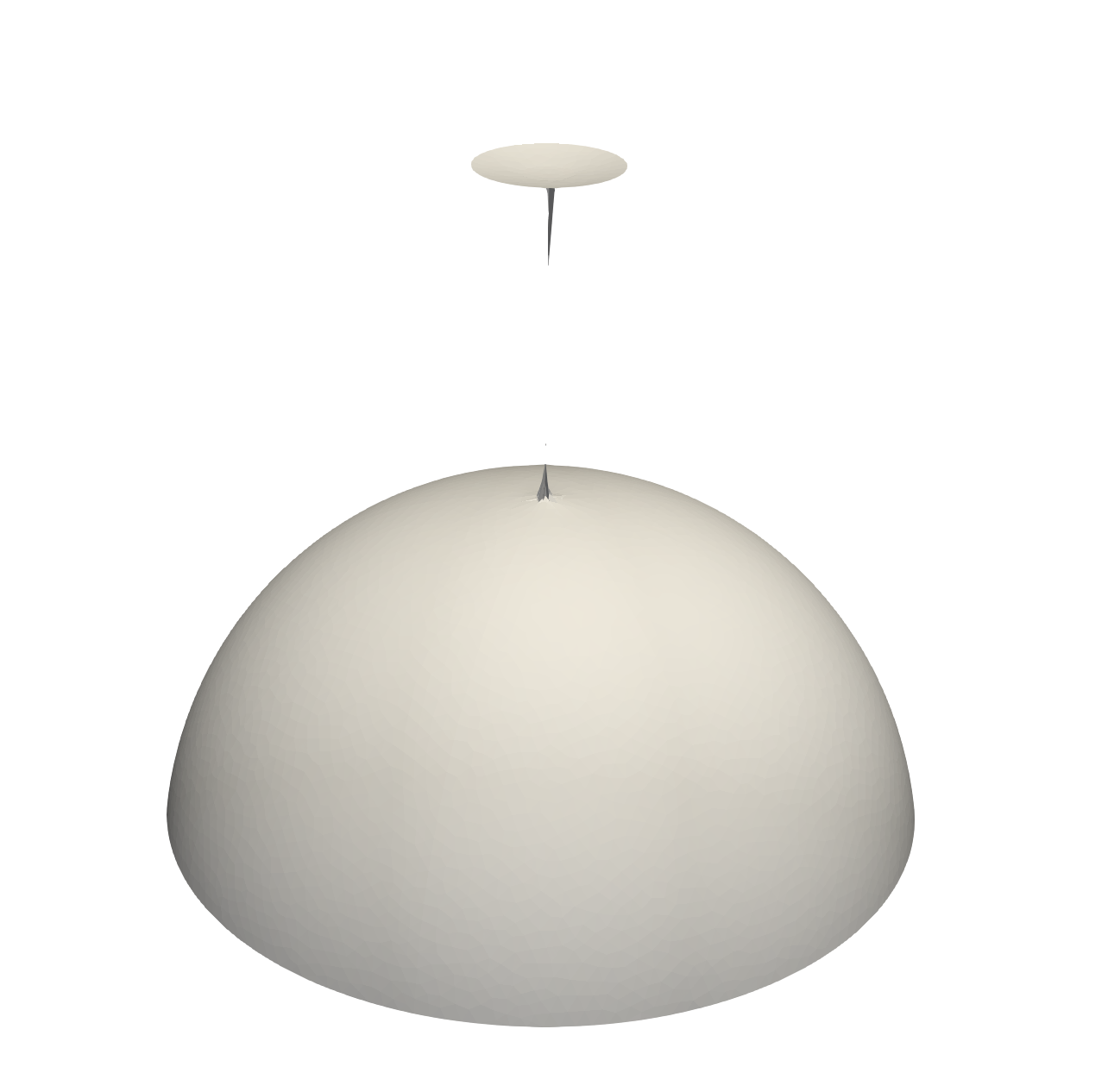}
\includegraphics[width=0.32\textwidth]{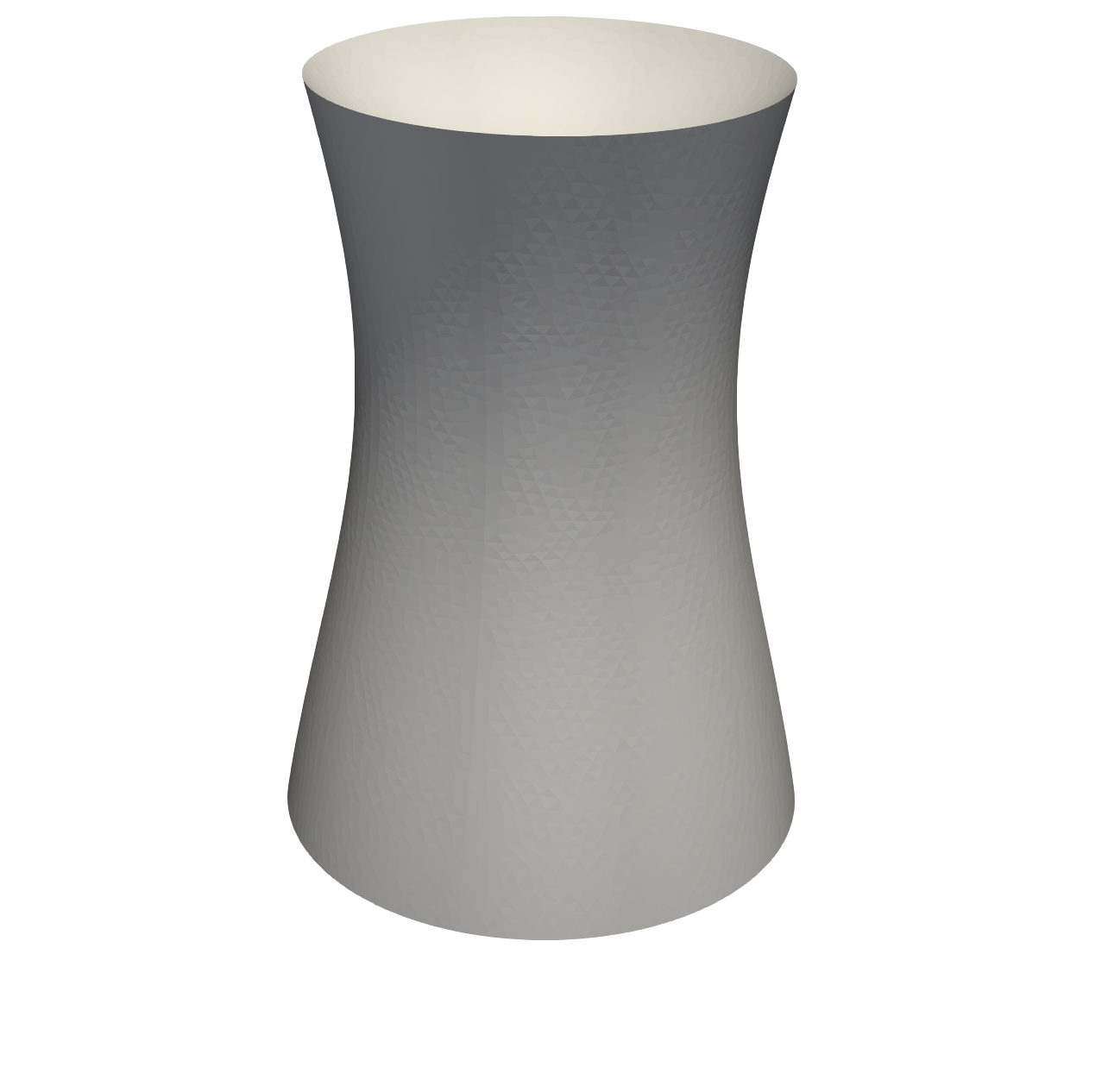}
\end{center}
\caption{Initial shape and evolution for the unstable saddle point from~\cite{rubinstein_2014} (left and middle). Proper local minimum on the right.}\label{fig:Rubinstein}
\end{figure}
Due to floating-point inaccuracies, a saddle point cannot be expected to be a stable point for a numerical gradient descent scheme. Indeed, the unduloid shape for a meniscus angle of $32^\circ$ proved to be almost stationary initially. However, a rapid descent of the objective was eventually achieved with the geometry indicating the desire to split the single capillary bridge into two droplets. This ultimately terminates the program as topological changes are not possible, yet. As expected, the Newton-iteration was unusable for this geometry due to the ill-conditioning of the Hessian. Contrary to this saddle-point, a Newton-scheme is possible for the initial geometry of the unduloid of meniscus angle $47.2^\circ$, where the provided starting guess of the discretized analytical shape quickly converges to a residuum of $7.940805\cdot 10^{-9}$, numerically identifying this shape as a proper local minimum. The resulting forces on the sphere were computed to be $F_p = -2.140941$ and $F_s = 4.481142$. The respective shapes are shown in Figure~\ref{fig:Rubinstein}.

\section{Summary and Conclusion}
A novel shape Newton approach for computing capillary bridges has been studied. In particular, we derive a variational formulation of the first and second order derivative of the Lagrangian of the total energy of liquid bridges between solid particles with respect to shape perturbations of the liquid phase. Special attention was paid on variational formulations, which respect the lower regularity of triangulated surfaces. Hence, expressions involving curvature are avoided. Furthermore, the desire to compute all quantities on shell meshes lead to novel expressions for the shape Hessian and the second order material derivative of the normal, which are also of interest for more general problems. The resulting expressions where confirmed to provide the same numerical values as the ``discretize-then-optimize'' approach, even for coarse meshes. Several test cases for different capillary bridges from the literature were revisited, confirming both the accuracy of the method as well as the performance of the Newton scheme.


\printbibliography

\end{document}